\documentclass{article}
\usepackage{spconf}
\usepackage{amsmath}
\usepackage{amssymb}
\usepackage{bm}
\usepackage{amsthm}
\usepackage{graphicx}
\usepackage{tikz}
\usetikzlibrary{arrows, automata}
\usepackage{pgfplots}
\usepgfplotslibrary{fillbetween}
\pgfplotsset{compat=1.5}
\pgfplotsset{major grid style={dashed}}
\usepackage[sorting=nty, style=ieee, citestyle=ieee]{biblatex}
\usepackage{algorithm}
\usepackage{algpseudocode}
\usepackage{algorithmicx}

\definecolor{Saffron}{HTML}{34282C}
\newtheorem{claim}{Claim}
\newtheorem{definition}{Definition}

\title{TROPICAL MODELING OF WEIGHTED TRANSDUCER ALGORITHMS ON GRAPHS}
\name{Emmanouil Theodosis and Petros Maragos}
\address{School of ECE, National Technical University of Athens, Athens, Greece\\
manostheodosis@gmail.com, maragos@cs.ntua.gr}
\begin{document}
\ninept
\maketitle
\begin{abstract}
Weighted Finite State Transducers (WFSTs) are versatile data structures that can model a great number of problems, ranging from Automatic Speech Recognition to DNA sequencing. Traditional computer science algorithms are employed when working with these structures in order to optimise their size, but also the runtime of decoding algorithms. However, these algorithms are not unified under a common framework that would allow for their treatment as a whole. Moreover, the inherent geometrical representation of WFSTs, coupled with the topology-preserving algorithms that operate on them make the structures ideal for tropical analysis. The benefits of such analysis have a twofold nature; first, matrix operations offer a connection to nonlinear vector space and spectral theory, and, second, tropical algebra offers a connection to tropical geometry. In this work we model some of the most frequently used algorithms in WFSTs by using tropical algebra; this provides a theoretical unification and allows us to also analyse aspects of their tropical geometry. Further, we provide insights via numerical examples.
\end{abstract}
\section{Introduction}
\label{sec:intro}
Weighted Finite State Transducers (WFSTs) are complex mathematical objects which find application in many fields ranging from language and speech processing to computational biology. There exists a multitude of algorithms that operate on WFSTs (\cite{Mohr00}, \cite{Mohr09}, \cite{MPR02}). The most prominent and most studied is the Viterbi algorithm (\cite{Forn73}, \cite{Rabi89}, \cite{Vite67}), which stems from the field of telecommunications. Others include the weight pushing and the epsilon removal algorithms, stemning from computer science. However, such algorithms are usually computationally expensive, which is undesirable for practical applications. Moreover, besides algorithms that simply utilise the WFST, there are more intrusive algorithms that alter its parameters in an effort to optimise subsequent decoding, while maintaining its inherent structure. Some of these algorithms aim to directly reduce the size of the states and arcs in WFST, and thus immediately affecting the time requirements of the decoding. On the other hand, certain algorithms try to indirectly affect the execution speed, by adapting the weights between states so that pruning algorithms examine fewer paths.

These algorithms admit modeling through tropical algebra and tropical geometry (\cite{Kuo06}, \cite{PaSt04b}, \cite{Mohr00}, \cite{Mohr09}, \cite{MPR02}), however no efforts have been made to thoroughly explore their tropical aspects beyond the expression of scalar arithmetic with operations from the tropical semiring. For detailed background on tropical algebra and the tropical semiring we refer the reader to \cite{Cuni79}, \cite{Butk10}, \cite{GoMi08}, \cite{Mara17}, \cite{Simo94}, and \cite{Pin98}. In this paper we model the algorithms using tropical algebra and matrix operations, resulting in novel expressions in closed matrix form. We also explain aspects of the geometry of certain algorithms, namely the Viterbi pruning.

\cite{Pin98} first introduces the min-plus arithmetic. References \cite{Mohr00}, \cite{Mohr09}, and \cite{MPR02} are some of the most influential of the field, studying the WFST structures and proposing the corresponding algorithms. In \cite{Vite67} the Viterbi algorithm was first introduced as an optimal decoding algorithm. Reference \cite{Cuni79} is a thorough study of nonlinear algebras, namely the minimax algebra. In \cite{Butk10} the author focuses on max-plus algebra from a control theory viewpoint. Max-plus algebra is also studied, along with its applications, in \cite{BCOQ93} and \cite{Gaub97}. In \cite{Mara17} the author offers a comprehensive study of systems on weighted lattices as a unification of max-plus algebra and its generalisations. References \cite{ChMa17} and \cite{ThMa18} from our group are efforts to model perceptrons in max-plus algebra, in the case of the former, and the Viterbi algorithm and its pruning variant in min-plus algebra, in the case of the latter.

In this paper we provide a theoretical unification of WFSTs algorithms by modeling them using tropical algebra, which also allows for their further analysis using tools from minimax matrix theory (\cite{Cuni79}). We first model the weight pushing algorithm, an non-intrusive algorithm that aims to speed up pruning by propagating the weights to earlier states of the WFST. Then we model the epsilon removal algorithm, which alters the structure of the WFST in order to remove unecessary states and trastitions, thus reducing its size and immediately affecting decoding. We present previous results regarding the modeling of the Viterbi algorithm and its pruning variant. Finally, we further explore the properties of certain metrics defined through the Viterbi pruning and elaborate on their motivation. Our modeling aspires to offer a connection with and unification via nonlinear vector space theory of weighted lattices (\cite{Mara17}) and aspires to allow for spectral analysis of these algorithms. In addition, we provide links with tropical geometry, similar to the efforts in \cite{ChMa17} and \cite{ThMa18}.

In Section \ref{sec:back} we present elements of tropical algebra that will be useful in our analysis. Section \ref{sec:mod} contains the modeling of tha various algorithms in tropical algebra; namely the weight pushing, epsilon removal, and Viterbi algorithms. Finally, in Section \ref{sec:anal} we revisit the geometry of the Viterbi pruning and we better explain the motivation for and the properties of metrics defined in previous work (\cite{ThMa18}).

\section{Background}
\label{sec:back}
Tropical algebra is similar to linear algebra. Like linear algebra studies systems of linear equations and their properties, tropical algebra studies systems of nonlinear equations (namely, min-plus equations) and their properties. Its main pair of operations is the pair $(\min, +)$, and we will use $\wedge$ to denote the minimum. The vectors and matrices of tropical algebra exist on the extended real multidimensional space defined by $\mathbb{R}_{\min} = \mathbb{R} \cup \{ +\infty\}$. In this paper, we follow the notation of \cite{Mara17} for the operations on weighted lattices. Let $\mathbf{A}, \mathbf{B} \in \mathbb{R}_{\min}^{n \times m}$. Then the min-plus product between these matrices, denoted by $\boxplus$, is given by:
\begin{equation}
    \label{eq:minplusmult}
    \left( \mathbf{A} \boxplus \mathbf{B}\right)_{ij} = \bigwedge_{k = 1}^{m} A_{ik} + B_{kj}
\end{equation}

We will also make extensive use of two very important matrices for tropical algebra. In partiular, we will use:
\begin{itemize}
    \item the matrix $\Gamma(\mathbf{A})$ of a matrix $\mathbf{A}$, defined as:
    \begin{equation}
        \label{eq:gamma}
        \Gamma(\mathbf{A}) = \mathbf{A} \wedge \mathbf{A}^2 \wedge ... \wedge \mathbf{A}^n \wedge ...
    \end{equation}
    \item the matrix $\Delta(\mathbf{A})$ of a matrix $\mathbf{A}$, defined as:
    \begin{equation}
        \label{eq:delta}
        \Delta(\mathbf{A}) = \mathbf{I} \wedge \mathbf{A} \wedge \mathbf{A}^2 \wedge ... \wedge \mathbf{A}^n \wedge ...
    \end{equation}
\end{itemize}

We can see that $\Delta(\mathbf{A}) = \mathbf{I} \wedge \Gamma(\mathbf{A})$. These two matrices are very important in tropical algebra, because they provide solutions to the eigenvector problems. In particular:
\begin{itemize}
    \item the matrix $\Gamma(\mathbf{A})$ provides solutions to the min-plus eigenvector-eigenvalue problem $\mathbf{A} \boxplus \mathbf{x} = \lambda \boxplus \mathbf{x}$.
    \item the matrix $\Delta(\mathbf{A})$ provides solutions to the generalised min-plus eigenvector-eigenvalue problem $\mathbf{A} \boxplus \mathbf{x} \geq \lambda \boxplus \mathbf{x}$.
\end{itemize}

Tropical geometry (\cite{Zieg12}, \cite{MaSt15}) aims to generalise the ideas of Euclidean geometry to the tropical setting. This proves useful in many cases because tropical curves are piecewise linear, which offers immediate bounds for the solution space of problems, but also offers ties to linear programming and its algorithms. Similar to its Euclidean counterpart, a tropical line is given by Equation \eqref{eq:tropline}:
\begin{equation}
    \label{eq:tropline}
    y = \alpha + x \wedge \beta = \min(\alpha + x, \beta)
\end{equation}
Similarly to the tropical lines we can define \emph{tropical halfspaces} as:
\begin{definition}
    Let $\mathbf{a}, \mathbf{b} \in \mathbb{R}_{\min}^{n+1}$. An \emph{affine tropical halfspace} is a subset of $\mathbb{R}_{\min}^{n}$ defined by:
    \begin{align*}
        T (\mathbf{a}, \mathbf{b}) := \{ \mathbf{x} \in \mathbb{R}_{\min}^{n}: && \left(\bigwedge_{i = 1}^{n} a_i + x_i \right) \wedge a_{n+1} \geq \\
        && \left(\bigwedge_{i = 1}^{n} b_i + x_i \right) \wedge b_{n+1} \}
    \end{align*}
\end{definition}
In the text we will reference tropical polytopes. These mathematical objects arise from the combination of tropical halfspaces:
\begin{definition}
    A bounded intersection of a finite number of tropical halfspaces is will be called a \emph{tropical polytope}.
\end{definition}

\section{Modeling}
\label{sec:mod}
\subsection{Weight pushing}
\label{subsec:wp}
\begin{figure}[t]
    \begin{center}
        \begin{tikzpicture}
            \node[initial, state] (0) at (-1, 0) {$0$};
            \node[state] (1) at (1.5, 1) {$1$};
            \node[state] (2) at (1.5, -1) {$2$};
            \node[accepting, state] (3) at (4, 1) {$3$};
            \node[state] (4) at (4, -1) {$4$};

            \draw[->] (0) edge node [above left] {$\alpha:A/1$} (1)
                    edge node [below left] {$\alpha:A/2$}  (2)
                (1) edge node [above] {$\zeta:Z/42$}  (3)
                (2) edge node [below] {$\chi:X/3$}  (4);
            \draw (4) -- (5, -0.5);
            \draw (4) -- (5, -1.5);
            \node (5) at (5.25, -0.5) {$...$};
            \node (6) at (5.25, -1.5) {$...$};

            \draw (1.9, -1.8) -- (1.9, -2.6);
            \draw (2.1, -1.8) -- (2.1, -2.6);
            \draw (1.8, -2.4) -- (2, -2.8) -- (2.2, -2.4);

            \node[initial, state] (7) at (-1, -4.5) {$0$};
            \node[state] (8) at (1.5, -3.5) {$1$};
            \node[state] (9) at (1.5, -5.5) {$2$};
            \node[accepting, state] (10) at (4, -3.5) {$3$};
            \node[state] (11) at (4, -5.5) {$4$};

            \draw[->] (7) edge node [above left] {$\alpha:A/43$} (8)
                    edge node [below left] {$\alpha:A/5$}  (9)
                (8) edge node [above] {$\zeta:Z/0$}  (10)
                (9) edge node [below] {$\chi:X/0$}  (11);
            \draw (11) -- (5, -5);
            \draw (11) -- (5, -6);
            \node (12) at (5.25, -5) {$...$};
            \node (13) at (5.25, -6) {$...$};
        \end{tikzpicture}
        \caption{WFST transducing a sequence of lowercase greek letters to the corresponding sequence of capital letters. The structure of the network allows for weights to be pushed to earlier transitions without altering the ranking and cost of the accepted paths.}
        \label{fig:wp}
    \end{center}
\end{figure}
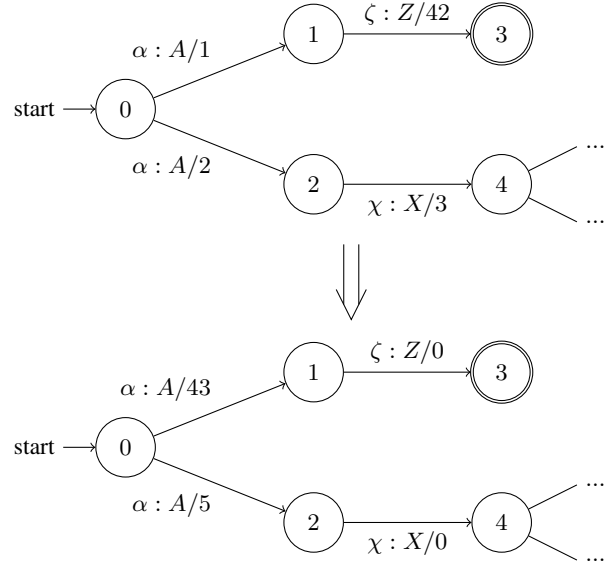

The weight pushing algorithm is an essential algorithm for practical application of the WFST framework. The algorithm aims to propagate the weights to earlier states of the structure, to the effect that low-probability paths are recognised earlier in the decoding sequence, and thus have a higher chance of being pruned by pruning algorithms. An irrevocable requirement is that the underlying structure of the WFST must remain the same: the algorithm might alter the weights, but the set of accepted paths and their total weights must stand unaffected. An example highlighting the weight pushing operation appears in Figure \ref{fig:wp}. An improbable path that has, at an early stage, a low cost will consume computational resources, where that could have been avoided by pushing the overall weight in earlier transitions.

The algorithm can be divided into two parts; a first part, where a potential (meaning the amount that can be propagated to earlier states) is calculated, and a second part where the actual update of the parameters occurs. A single iteration of the traditional algorithm for calculating the potential can be written in the form:
\begin{equation}
    \label{eq:pot}
    \mathbf{v}_{i + 1} = \mathbf{v}_{i} \wedge \mathbf{A} \boxplus \mathbf{v}_{i}
\end{equation}
where $\mathbf{v}_{i}$ is the potential vector for the $i$-th iteration, and $\mathbf{v}_{0} = \pmb{\rho}$, where $\pmb{\rho}$ is the emission vector. By recursively substituting the values, we get that the final value of the potential vector is:
\begin{equation}
    \label{eq:pot_final}
    \mathbf{v}_{N} = \pmb{\rho} \wedge \mathbf{A} \boxplus \pmb{\rho} \wedge \mathbf{A}^2 \boxplus \pmb{\rho} \wedge ... \wedge \mathbf{A}^n \boxplus \pmb{\rho} \wedge ... = \Delta(\mathbf{A}) \boxplus \pmb{\rho}
\end{equation}

\begin{claim}
    The calculation of Equation \eqref{eq:pot_final}
    \begin{equation}
        \mathbf{v}_N = \Delta(\mathbf{A}) \boxplus \pmb{\rho}
    \end{equation}
    is finite and $\Delta(\mathbf{A}) = \mathbf{I} \wedge \mathbf{A} \wedge ... \wedge \mathbf{A}^{n-1}$.
\end{claim}

The claim is proven by the fact that we have assumed that there aren't any cycles of negative length in the WFST, and such the shortest paths between every pair of states are finite.

Having computed the potential vectors, we define four diagonal matrices that will be useful for updating the parameters of the WFST. In particular:
\begin{itemize}
    \item The matrix $\pmb{\Lambda}$ of the input weights, whose diagonal is the input weight vector $\pmb{\lambda}$.
    \item The matrix $\mathbf{V}_N$ of the potentials, whose diagonal is the potential vector $\mathbf{v}_N$.
    \item The matrix $\mathbf{V}_{-N}$ of the negative potentials, whose diagonal is the negative of the potential vector $\mathbf{v}_N$.
    \item The matrix $\mathbf{P}$ of the emission weights, whose diagonal is the emission weight vector $\pmb{\rho}$.
\end{itemize}

Having defined these matrices, the updated parameters of the WFST are as follows:
\begin{equation}
    \pmb{\lambda}' = \pmb{\Lambda} \boxplus \mathbf{v}_N, \quad \pmb{\rho}' = \mathbf{P} \boxplus \mathbf{v}_{-N}, \quad \mathbf{A}' = \mathbf{V}_{-N} \boxplus \mathbf{A} \boxplus \mathbf{V}_{N}
\end{equation}

\subsection{Epsilon removal}
\label{subsec:eps}
Epsilon removal is an algorithm that aims to reduce the number of states and transitions in the WFST, while maintaining its underlying structure, in order to reduce the running time of the Viterbi algorithm. To accomplish that, an effort is made to reduce epsilon transitions (meaning transtitions with no input or output symbols).

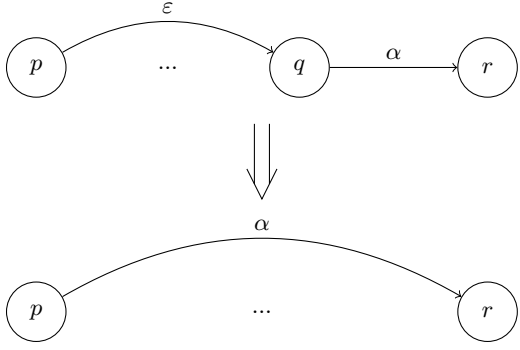
\begin{figure}[t]
    \begin{center}
        \begin{tikzpicture}
            \node[state] (0) at (-1, -0.25) {$p$};
            \node[state] (1) at (2.5, -0.25) {$q$};
            \node[state] (2) at (5, -0.25) {$r$};

            \draw[->] (0) edge [bend left] node [above] {$\varepsilon$} (1)
                (1) edge node [above] {$\alpha$} (2);
            \node (3) at (0.75, -0.25) {$...$};

            \draw (1.9, -1) -- (1.9, -1.8);
            \draw (2.1, -1) -- (2.1, -1.8);
            \draw (1.8, -1.6) -- (2, -2) -- (2.2, -1.6);

            \node[state] (4) at (-1, -3.5) {$p$};
            \node[state] (5) at (5, -3.5) {$r$};

            \draw[->] (4) edge [bend left] node [above] {$\alpha$} (5);
            \node (7) at (2, -3.5) {$...$};
        \end{tikzpicture}
        \caption{Epsilon transitions increase the total number of states and transitions of the WFST, increasing its size and reducing its efficiency. Removing states and replacing transitions, while maintaining the same accepted paths and costs, improves the runtime of decoding algorithms.}
        \label{fig:epsilon_ex}
    \end{center}
\end{figure}

The traditional algorithm for epsilon removal can be illustrated in Figure \ref{fig:epsilon_ex}. In essence, the algorithm computes, for each state $p$ its epsilon closure, and then adds transitions from the state $p$ to each state reachable from states in the epsilon closure.

To model the traditional epsilon removal algorithm in tropical algebra we need to define two matrices, in addition to the transition matrix $\mathbf{A}$ of the WFST:
\begin{itemize}
    \item the matrix $\mathbf{\Sigma}_I$, which is the input symbol matrix and $\left(\Sigma_I\right)_{ij}$ contains the input symbol for the transition of the state $i$ to state $j$.
    \item the matrix $\mathbf{\Sigma}_O$, which is the output symbol matrix and $\left(\Sigma_O\right)_{ij}$ contains the output symbol for the transition of the state $i$ to state $j$.
\end{itemize}

For completeness sake, we need to make two remarks before we proceed to the modeling:
\begin{itemize}
    \item We only consider as epsilon transitions ones where both the input and the output symbols are $\varepsilon$. This is a very common assumption in the field, and usually a synchronization algorithm has already been performed, in order to better match input and output $\varepsilon$.
    \item We assume that there can only be a single transition between two states, regardless of whether there exist transitions with different symbols or weights. While this might seem restrictive, in practice it isn't, and can even be circumvented.
\end{itemize}

We need to define another two matrices in order to model epsilon removal, which make up the transition matrix $\mathbf{A}$:
\begin{align}
    \label{eq:Ematrix}
    E_{ij} &=
        \begin{cases}
            \alpha_{ij}, & \text{if } \left(\Sigma_I\right)_{ij} = \left(\Sigma_O\right)_{ij} = \varepsilon \\
            \infty,   & \text{otherwise}
        \end{cases},
        \\
    \left(A_{\varepsilon}\right)_{ij} &=
        \begin{cases}
            \alpha_{ij}, & \text{if } \left(\Sigma_I\right)_{ij} = \left(\Sigma_O\right)_{ij} \neq \varepsilon \\
            \infty,   & \text{otherwise}
        \end{cases}
\end{align}

Essentially, we decompose matrix $\mathbf{A}$ using tha matrices of Equation \eqref{eq:Ematrix}. We can see that $\mathbf{A} =\mathbf{A}_{\varepsilon} \wedge \mathbf{E}$.

\begin{claim}
    Let $\mathbf{E}$ the matrix defined in \eqref{eq:Ematrix}. Then, the matrix
    \begin{equation}
        \Gamma(\mathbf{E}) = \mathbf{E} \wedge \mathbf{E}^2 \wedge ... \wedge \mathbf{E} ^ n \wedge ...
    \end{equation}
    is finite and equal to $\Gamma(\mathbf{E}) = \mathbf{E} \wedge \mathbf{E}^2 \wedge ... \wedge \mathbf{E}^{n-1}$, and moreover expresses the epsilon closure for all the states of the WFST.
\end{claim}

The claim is proven by the fact that an inherent assumption in WFSTs is that there aren't any cycles of negative weight (and thus the shortest distances are finite). Since there aren't any cycles of negative weight in the original WFST, there aren't any such cycles in the WFST where we kept only the epsilon transitions. Having the epsilon closure of each state, the updated transition matrix and emission vector are simply the tropical addition (that is, the minimum) between the previous values and the values that emerge from the epsilon closure. In particular, the new transition matrix $\mathbf{A}'$ takes the form:
\begin{equation}
    \mathbf{A}' = \mathbf{A}_\varepsilon \wedge \left(\Gamma(\mathbf{E}) \boxplus \mathbf{A}_\varepsilon \right) = \Delta(\mathbf{E}) \boxplus \mathbf{A}_\varepsilon
\end{equation}
whereas the new emission vector $\pmb{\rho}'$ takes the form:
\begin{equation}
    \pmb{\rho}' = \pmb{\rho} \wedge \left(\Gamma(\mathbf{E}) \boxplus \pmb{\rho}\right) = \Delta(\mathbf{E}) \boxplus \pmb{\rho}
\end{equation}

\subsection{Viterbi}
The Viterbi algorithm aims to decode a sequence of input symbols, meaning that it tries to map the sequence of symbols to the sequence of states that has the highest probability. Formally, it is known that the Viterbi algorirthm can be written in the following max-product form:
\begin{equation}
    \label{eq:viterbi}
    q_i(t) = \left( \max_{j} w_{ji} q_j(t-1) \right) \cdot b_i (\sigma_t)
\end{equation}
where $w_{ji}$ is the probability of transitioning from state $j$ to state $i$, $b_i(\sigma_t)$ denotes the observation probability of the symbol $\sigma_t$ at state $i$, and, finally, $q_i(t)$ is the maximum probability for that current state, calculated along the path from the previous states. In \cite{ThMa18} we postulated that the Viterbi algorithm can be written in a closed matrix form in tropical algebra as:
\begin{equation}
    \label{eq:trop_viterbi}
    \mathbf{x}(t) = \mathbf{P}(\sigma_t) \boxplus \mathbf{A}^T \boxplus \mathbf{x}(t-1)
\end{equation}
where $\mathbf{x}(t) = -\log\mathbf{q}(t)$, $\mathbf{A} = -\log\mathbf{W}$, and $\mathbf{P}(\sigma_t)$ is a diagonal matrix whose diagonal is the vector $\mathbf{p}(\sigma_t)$, with $\mathbf{p}(\sigma_t) = -\log\mathbf{b}(\sigma_t)$.

\subsection{Viterbi pruning}
\begin{figure}[t]
    \begin{center}
        \begin{tikzpicture}
            \draw[thick, opacity=0.25] (-1.4, 2.75) circle (10pt);
            \draw[thick, opacity=0.25] (-1.4, 2) circle (10pt);
            \draw[thick, opacity=0.25] (-1.4, 1.25) circle (10pt);

            \draw[thick] (0.3, 2.75) circle (10pt);
            \draw[thick] (0.3, 2) circle (10pt);
            \draw[thick] (0.3, 1.25) circle (10pt);

            \draw[thick, opacity=0.25] (2, 2.75) circle (10pt);
            \draw[thick, opacity=0.25] (2, 2) circle (10pt);
            \draw[thick, opacity=0.25] (2, 1.25) circle (10pt);

            \draw[->, opacity=0.25] (-1.05, 2.75) -- (-0.05, 2.75);
            \draw[->, opacity=0.25] (-1.05, 2.75) -- (-0.05, 2);
            \draw[->, opacity=0.25] (-1.05, 2.75) -- (-0.05, 1.25);
            \draw[->, opacity=0.25] (-2, 2.75) -- (-1.75, 2.75);
            \draw[->, opacity=0.25] (-2, 2.55) -- (-1.75, 2.75);
            \draw[->, opacity=0.25] (-2, 2.35) -- (-1.75, 2.75);

            \draw[->, opacity=0.25] (-1.05, 2) -- (-0.05, 2.75);
            \draw[->, opacity=0.25] (-1.05, 2) -- (-0.05, 2);
            \draw[->, opacity=0.25] (-1.05, 2) -- (-0.05, 1.25);
            \draw[->, opacity=0.25] (-2, 2.2) -- (-1.75, 2);
            \draw[->, opacity=0.25] (-2, 2) -- (-1.75, 2);
            \draw[->, opacity=0.25] (-2, 1.8) -- (-1.75, 2);

            \draw[->, opacity=0.25] (-1.05, 1.25) -- (-0.05, 2.75);
            \draw[->, opacity=0.25] (-1.05, 1.25) -- (-0.05, 2);
            \draw[->, opacity=0.25] (-1.05, 1.25) -- (-0.05, 1.25);
            \draw[->, opacity=0.25] (-2, 1.65) -- (-1.75, 1.25);
            \draw[->, opacity=0.25] (-2, 1.45) -- (-1.75, 1.25);
            \draw[->, opacity=0.25] (-2, 1.25) -- (-1.75, 1.25);

            \draw[->] (0.65, 2.75) -- (1.65, 2.75);
            \draw[->] (0.65, 2.75) -- (1.65, 2);
            \draw[->] (0.65, 2.75) -- (1.65, 1.25);
            \draw[opacity=0.25] (2.35, 2.75) -- (2.6, 2.75);
            \draw[opacity=0.25] (2.35, 2.75) -- (2.6, 2.55);
            \draw[opacity=0.25] (2.35, 2.75) -- (2.6, 2.35);

            \draw[->] (0.65, 2) -- (1.65, 2.75);
            \draw[->] (0.65, 2) -- (1.65, 2);
            \draw[->] (0.65, 2) -- (1.65, 1.25);
            \draw[opacity=0.25] (2.35, 2) -- (2.6, 2.2);
            \draw[opacity=0.25] (2.35, 2) -- (2.6, 2);
            \draw[opacity=0.25] (2.35, 2) -- (2.6, 1.8);

            \draw[->] (0.65, 1.25) -- (1.65, 2.75);
            \draw[->] (0.65, 1.25) -- (1.65, 2);
            \draw[->] (0.65, 1.25) -- (1.65, 1.25);
            \draw[opacity=0.25] (2.35, 1.25) -- (2.6, 1.65);
            \draw[opacity=0.25] (2.35, 1.25) -- (2.6, 1.45);
            \draw[opacity=0.25] (2.35, 1.25) -- (2.6, 1.25);

            \draw (-0.15, 3.2) -- (0.75, 3.2) -- (0.75, 0.8) -- (-0.15, 0.8) -- cycle;
            \draw (-0.15, 0.8) -- (0.3, 0.25) -- (0.75, 0.8);

            \draw[thick] (0, 0) -- (-0.05, 0) -- (-0.05, -2.05) -- (0, -2.05);
            \draw[thick] (0.6, 0) -- (0.65, 0) -- (0.65, -2.05) -- (0.6, -2.05);

            \node (0) at (0.3, -0.25) {$0.57$};
            \draw[red, opacity = 0.8] (0.3, -0.25) circle (9pt);

            \node (1) at (0.3, -1.75) {$0.72$};
            \draw[blue, opacity=0.8] (0.3, -1.75) circle (9pt);

            \node (2) at (0.3, -1) {$0.86$};
            \draw[green, opacity=0.8] (0.3, -1) circle (9pt);

            \draw[->, red, opacity=0.8] (0.75, -0.25) to [bend left] (2.8, -0.5);
            \draw[->, green, opacity=0.8] (0.75, -1) to [bend left] (3.25, -3.9);
            \draw[->, blue, opacity=0.8] (0.3, -2.1) to [bend right] (2.25, -4.5);

            \draw[red, opacity=0.8] (3, 0) -- (3, -4);
            \draw[blue, opacity=0.8] (3, -4) -- (2, -5);
            \draw[green, opacity=0.8] (3, -4) -- (4.5, -4);

            \draw[dashed, opacity=0.5] (2, -5) -- (3.5, -5);
            \draw[dashed, opacity=0.5] (3, 0) -- (4.5, 0);
            \draw[dashed, opacity=0.5] (2, -1) -- (3.5, -1);

            \draw[dashed, opacity=0.5] (4.5, -4) -- (3.5, -5);
            \draw[dashed, opacity=0.5] (3, 0) -- (2, -1);
            \draw[dashed, opacity=0.5] (4.5, 0) -- (3.5, -1);

            \draw[dashed, opacity=0.5] (4.5, 0) -- (4.5, -4);
            \draw[dashed, opacity=0.5] (2, -1) -- (2, -5);
            \draw[dashed, opacity=0.5] (3.5, -1) -- (3.5, -5);

            \node[text=red] (3) at (2.75, -1.5) {$r_1$};
            \node[text=blue] (4) at (2.8, -4.5) {$r_3$};
            \node[text=green] (5) at (4, -3.75) {$r_2$};

        \end{tikzpicture}
        \caption{At each step along the trellis, the state vector $\mathbf{x}(t)$ and the leniency vector $\pmb{\eta}$ of Equation \eqref{eq:lower} define a polytope. The three values of the vector $\mathbf{r}$ are shown with colors (see Equation \eqref{eq:nu}).}
        \label{fig:vit_pr}
    \end{center}
\end{figure}
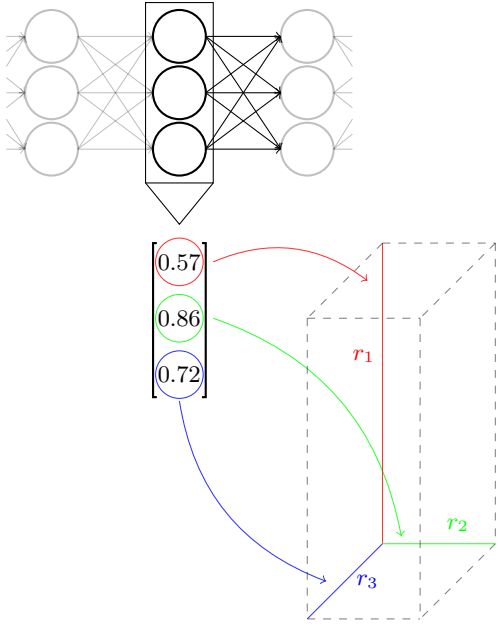
The Viterbi pruning is a variant of the Viterbi algorithm that aims to sacrifice the optimality of decoding in an effort to significantly speed up the decoding process. Usually, pruning is based on one of the following criteria, or even their combination:
\begin{itemize}
    \item users determine a leniency parameter $\theta$, and at each step only the paths that are at most $\theta$ from the optimal path survive.
    \item users determine a beam width $\kappa$, and at each step only the $\kappa$-best paths survive the pruning.
\end{itemize}

In \cite{ThMa18} we modeled the Viterbi pruning in tropical algebra using \emph{Cuninghame-Green's inverse} (\cite{Cuni79}). Therein it is proven that the negative elements of
\begin{equation}
    \label{eq:cun}
    \overline{\mathbf{y}} = \mathbf{X}^{\#}(t)\boxplus'\pmb{\eta}
\end{equation}
indicate the indices that need to be pruned. The matrix $\mathbf{X}(t)$ is a diagonal matrix whose diagonal is the state vector $\mathbf{x}(t)$, and $\mathbf{X}^{\#}(t):= -\mathbf{X}^T(t)$. Also, $\pmb{\eta} = \theta + \frac{1}{2}\left(\mathbf{x}^T(t) \boxplus \mathbf{x}(t) \right) + \mathbf{0}$, where $\theta$ is the leniency parameter and $\mathbf{0}$ is a vector that comprises solely of 0. Finally, $\boxplus'$ denotes the max-plus matrix multiplication.

Moreover, if we consider a variable vector $\mathbf{z}$ and bound it using:
\begin{itemize}
    \item the Viterbi update law of Equation \eqref{eq:viterbi}, and thus:
    \begin{equation}
        \label{eq:upper}
        \mathbf{z} \geq \mathbf{b}, \quad \mathbf{b} = \mathbf{P}(\sigma_t)\boxplus\mathbf{A}^T\boxplus\mathbf{x}(t-1)
    \end{equation}
    \item the pruning vector of Equation \eqref{eq:cun}, and thus
    \begin{equation}
        \label{eq:lower}
        \mathbf{z} \leq \pmb{\eta}, \quad \pmb{\eta} = \theta + \frac{1}{2}\left(\mathbf{b}^T\boxplus\mathbf{b}\right) + \mathbf{0}
    \end{equation}
\end{itemize}
Then, the combination of Equations \eqref{eq:upper} and \eqref{eq:lower} defines a tropical polytope for each step of the Viterbi algorithm. Then, for each iteration two metrics are defined based on that polytope:
\begin{itemize}
    \item a normalised volume metric $\nu$:
    \begin{equation}
        \label{eq:nu}
        \nu = -\frac{1}{\mathrm{supp}(\mathbf{z})} \sum_{i \in \mathrm{supp}(\mathbf{z})} \frac{\log r_i}{\log \left( \max \mathbf{r}\right)}
    \end{equation}
    \item a normalised entropy metric $\varepsilon$:
    \begin{equation}
        \label{eq:vare}
        \varepsilon = -\frac{1}{\mathrm{supp}(\mathbf{z})} \sum_{i \in \mathrm{supp}(\mathbf{z})}{-z_i(t) \cdot e^{-z_i(t)}}
    \end{equation}
\end{itemize}
where $r_i = \eta - z_i$. Essentially, $r_i$ is the degree to which each dimension satisfies the Viterbi constraints.

\section{Discussion of Geometry}
\label{sec:anal}
We devote this section to the further analysis of the metrics of Equations \eqref{eq:nu} and \eqref{eq:vare}, and also the motivation behind their definition. At every iteration of the Viterbi pruning algorithm consider the state vector $\mathbf{x}(t)$ along with the leniency vector $\pmb{\eta}$ of Equation \eqref{eq:lower}. In unison, these vectors define a tropical polytope for each iteration of the algorithm. The indices of the state vector that satisfy the constraints imposed by the leniency vector act as the sides of this polytope, and the difference between the value of the leniency vector and the state vector at that index constitute the vector $\mathbf{r}$ of Equation \eqref{eq:nu}. Figure \ref{fig:vit_pr} visualises the polytope of each iteration, and also highlights the vector $\mathbf{r}$. Duscussing the metrics further:
\begin{itemize}
    \item Consider the normalised volume of \eqref{eq:nu}. The metric $\nu$ can offer a quantitative estimate of the solution space that the Viterbi pruning admits. Indeed, since values of $r_i$'s in Equation \eqref{eq:nu} are normalised, utilising this metric can provide a measure of how many paths the current choice of the leniency parameter $\theta$ allows to survive. Exploiting that remark, it is possible to monitor how this metric evolves throughout iterations, and adapt, when needed, the value of the leniency paramater $\theta$ in order to maintain a desired level of normalised volume.
    \item Consider the normalised entropy of \eqref{eq:vare}. The metric $\varepsilon$ can offer a qualitative estimate of the solution space that the Viterbi pruning admits. In information theory, entropy expresses the current degree of surprise incured by the observation of a sample. In essence, if the sample abides by the existing modeling of the assumed distribution, then it will have low entropy, as its value is in an expected range. However, if the sample has a significantly different value than those expected by the assumptions for the distribution, then the sample will have very high entropy, indicating that there may be an error in the original modeling of the distribution.
\end{itemize}

Thus, by utulising the above metrics we aim to reason about the solution space of the Viterbi pruning in two ways; a quantitative analysis of the relative size of the solution space, and a qualitative analysis of the likelihood of the paths of the solution space. Having such measures, we can examine how the solution space, and the quantity/quality of these solutions evolves over the execution of the Viterbi algorithm. Even more, we can introduce them to the design of the algorithm, so that the leniency parameter $\theta$ gets adapted to the needs of each iteration.

\section{Conclusion}
\label{sec:conc}
In this work we modeled algorithms that operate on WFSTs using tropical algebra and matrix operations on weighted lattices, unifying them under a common framework. First, we modeled the weight pushing algorithm by expressing the potential calculation as an instrumental matrix of tropical algebra. We then proceeded to model the epsilon removal algorithm by exploiting the $\min$-superposition of tropical algebra and expressing the epsilon closure as another important matrix in tropical algebra. Finally, we analysed some geometrical aspects of the Viterbi pruning, elaborating on metrics that were defined on previous work. In future work we aim to explore the connection of the structures with nonlinear vector space theory of weighted lattices and nonlinear spectral theory.

\newpage


\printbibliography

\end{document}